\documentclass[centertags,12pt]{article}

\usepackage{amsmath}
\usepackage{amsthm}
\usepackage{amscd}
\usepackage{amssymb}
\usepackage{verbatim}

\title{Hidden symmetries and arithmetic manifolds}
\author{Benson Farb and Shmuel Weinberger \thanks{Both authors are
supported in part by the NSF.}\\
{\it \small Dedicated to the memory of Robert Brooks}}

\numberwithin{equation}{section}

\theoremstyle{plain}
\newtheorem{theorem}{Theorem}[section]
\newtheorem{proposition}[theorem]{Proposition}

\newtheorem{conjecture}[theorem]{Conjecture}
\newtheorem{question}[theorem]{Question}

\theoremstyle{definition}
\newtheorem{definition}[theorem]{Definition}
\newtheorem{remark}[theorem]{Remark}

\def\remark{{\bf {\medskip}{\noindent}Remark. }}

\def\remarks{{\bf {\medskip}{\noindent}Remarks. }}
\def\endproof{$\diamond$ \bigskip}

\def\title{\em}

\newcommand\R{\mbox{\bf R}}

\newcommand\Z{\mbox{\bf Z}}
\newcommand\Q{\mbox{\bf Q}}

\DeclareMathOperator{\Out}{Out}

\DeclareMathOperator{\Rank}{Rank}

\DeclareMathOperator{\Fix}{Fix}
\DeclareMathOperator{\SL}{SL}

\DeclareMathOperator{\Isom}{Isom}

\DeclareMathOperator\Sol{Sol}
\DeclareMathOperator\Comm{Comm}

\renewcommand\to{\longrightarrow}

\begin{document}
\maketitle

\section{Introduction}

Let $M$ be a closed, locally symmetric Riemannian manifold of
nonpositive curvature with no local torus factors; for example take $M$
to be a hyperbolic manifold.  Equivalently, $M=K\backslash G/\Gamma$
where $G$ is a semisimple Lie group and $\Gamma$ is a cocompact lattice 
in $G$. For simplicity, we will always assume that $\Gamma$ is {\em
irreducible}, or equivalently that $M$ is not finitely covered by a
smooth product; we will also assume for simplicity that $\dim(M)>2$.  

Let $g_{loc}$ denote a locally symmetric metric on $M$.  The
Mostow Rigidity Theorem states that $g_{loc}$ is
unique up to homothety of the direct factors of $\widetilde{M}$
\footnote{We will henceforth abuse notation and call $g_{loc}$
the ``unique locally symmetric metric on $M$''.}.  
As the following result (essentially due to Borel) shows, 
the metric $g_{loc}$ is also special because it
reflects every symmetry of every other Riemannian metric on $M$.
\enlargethispage{.4in}

\begin{proposition}[Most symmetry]
\label{proposition:biggest}
Let $M$ be a closed, irreducible, locally symmetric Riemannian
$n$-manifold, $n>2$, of
nonpositive curvature with no local torus factors.  Then 
for every Riemannian metric $h$ on $M$, the group 
$\Isom(M,h)$ is isomorphic to a subgroup of $\Isom(M,g_{loc})$.
\end{proposition}

\remarks
\begin{enumerate}
\item By an old result of Bochner, the group 
$\Isom(M,g_{loc})$ is finite.  

\item It is {\em not} true in general 
that $\Isom(M,h)$ is topologically conjugate to a subgroup of 
$\Isom(M,g_{loc})$: there are examples of hyperbolic $n$-manifolds,
$n>4$, with periodic diffeomorphisms not topologically conjugate into 
$\Isom(M,g_{loc})$; we give examples in \S\ref{subsection:exotic}. 
By averaging, one sees that 
these diffeomorphisms preserve some Riemannian metric.

\item We conjecture that Proposition \ref{proposition:biggest} 
remains true when $M$
has finite volume, and the metric $h$ is any finite volume Riemannian
metric.  However, we are only able to prove this when $M$ admits a
locally symmetric metric of negative curvature.  
\end{enumerate}

In \S\ref{section:3proofs} below we give three proofs of
Propsition \ref{proposition:biggest}: 
one using coarse geometry and Smith theory, one
using harmonic maps, and one (essentially due to Borel) 
using algebraic topology.  These proofs are
meant as a warm-up for the proof of our main theorem 
(Theorem \ref{theorem:characterize1} 
below); they are also meant 
as an introduction to some of the ideas in \cite{FW}, where we will 
make use of the interplay among all three proofs.  

We remark that both the harmonic maps proof and the
algebraic-topological proof apply {\it
mutatis mutandis} to any aspherical manifold with a degree one map to
$M$ that is an isomorphism on fundamental group.  We also remark that
each of the proofs applies to give a slightly modified version of 
Proposition \ref{proposition:biggest} in dimension $2$.  

\bigskip
\noindent
{\bf Hidden symmetries. } Does the property of $g_{loc}$ given in
Proposition \ref{proposition:biggest} characterize the locally symmetric metric
up to homothety? The answer is obviously no: since $\Isom(M,g_{loc})$
is finite, it is easy to equivariantly perturb the metric (say in the
neighborhood of a single free orbit, for simplicity) to give a
non-isometric manifold with the same isometry group.  Indeed, there is 
an infinite-dimensions worth of such metrics.

The main discovery of this paper is that we can actually detect
$g_{loc}$ exactly by looking at the isometry groups of finite-sheeted covers
of $M$. To make this more precise, we will need the following.

\begin{definition}[Hidden symmetry]
Let $M$ be a finite volume Riemannian manifold.  A {\em hidden symmetry}
of $M$ is an isometry $\phi$ of a $d$-sheeted, $1<d<\infty$, 
Riemannian cover of $M$
with the property that $\phi$ is not the lift of an isometry.
\end{definition}

So, for example, the deck transformations on any finite cover of $M$,
which are just lifts of the identity, are not hidden symmetries of $M$.

To find examples of locally symmetric manifolds with a lot of hidden
symmetry, we must consider {\em arithmetic manifolds}.  These are
manifolds $M=\Gamma\backslash G/K$ with $G$ semisimple, $K\subset G$ a
maximal compact, and $\Gamma$ an arithmetic group.  
Borel showed that these exist in 
abundance in every semisimple Lie group $G$, and 
Margulis's Arithmeticity Theorem (see,
e.g., \cite{Ma}) states that if 
$\Rank_{\R}(G)>1$ and if $\Gamma$ is irreducible, then every such $M$  
is arithmetic.

Recall that the {\em commensurator} of a lattice $\Gamma$ in a Lie group
$G$ is defined to be
$$\Comm_G(\Gamma):=\{g\in G: g\Gamma g^{-1}\cap \Gamma \mbox{\ has finite
index in both\ }g\Gamma g^{-1} \mbox{\ and \ }\Gamma\}$$

Borel and Harish-Chandra proved that arithmetic groups have infinite
index in their commensurators.  Using this we will prove the following.

\begin{proposition}[Arithmetic metrics have hidden symmetry]
\label{proposition:arithmetic}
Let $M$ be a closed, arithmetic manifold.  Then $(M,g_{loc})$ has 
infinitely many hidden symmetries; these can
be taken to occur on irregular covers.
\end{proposition}

We remark that, while every lift to $\widetilde{M}=G/K$ of a hidden symmetry
lies in the commensurator $\Comm_G(\Gamma)$, the set of hidden
symmetries is a rather small subset of
$\Comm_G(\Gamma)$; see \S\ref{subsection:arithmetic:hidden} below.  
Margulis proved that irreducible, 
arithmetic lattices in semisimple Lie groups are
precisely those irreducible lattices that have infinite 
index in their commensurator (see \cite{Ma}).  
The following theorem can be viewed as a
generalization of his result.

\begin{theorem}[Characterizing the arithmetic metric]
\label{theorem:characterize1}
Let $M$ be an arithmetic manifold.  Then $g_{loc}$ is the unique metric
(up to homothety) on $M$ having infinitely many hidden symmetries.
Conversely, if $M=\Gamma\backslash G/K$ with $G$ semisimple and 
$\Gamma$ an irreducible, cocompact {\em non-arithmetic} lattice, then $M$ has finitely many hidden
symmetries in every metric, and there are uncountably many
homothety classes of metrics on $M$ which have the maximal number of
hidden symmetries.
\end{theorem}

Theorem \ref{theorem:characterize1} characterizes the locally
symmetric metrric via a countable amount of data, namely the (finite)
isometry groups of the finite sheeted covers of $M$.  
In the special case when $(M,h)$ is nonpositively curved, then Theorem
\ref{theorem:characterize1} is a 
consequence of theorems of Eberlein [Eb] and Margulis (see 
below).  A different approach from ours to proving Theorem
\ref{theorem:characterize1} might be to use 
Furman's main theorem in \cite{Fu} classifying 
which locally-compact topological groups contain a lattice in a 
semisimple Lie group as a lattice.

In \cite{FW} we will prove the following much stronger theorem:

\begin{theorem}[Characterization of arithmetic manifolds]
\label{theorem:general}
Let $(M,h)$ be any closed, aspherical Riemannian manifold. Suppose that
$\pi_1(M)$ contains no normal infinite abelian subgroup, and that $M$
is ``irreducible'' in the sense that no finite cover of $M$ is a smooth
product.  If $M$ has an infinitely many hidden symmetries, then $(M,h)$
is isometric to an arithmetic manifold.
\end{theorem}

\remarks 

\begin{enumerate}

\item Each of the hypotheses in Theorem \ref{theorem:general} is
necessary: every irreducible lattice in a semisimple Lie group is not
virtually a product, and does not contain any normal infinite abelian
subgroup.

\item Note that 
Theorem \ref{theorem:characterize1} above requires the assumption
that the fundamental group of $M$ is arithmetic; we make no such a
priori assumption in Theorem \ref{theorem:general}.  

\item In \cite{FW} we
give a similar characterization of (not necessarily arithmetic) locally
symmetric manifolds in terms of symmetries that potentially remain
hidden until one gets to the universal cover.

\item The special role of normal abelian subgroups can be seen by
considering, for example, any $3$-dimensional solvmanifold $M$; this is
a manifold whose universal
cover is the unique $1$-connected non-nilpotent, solvable Lie group 
$\Sol$.  Every such $M$ has an infinite number of hidden symmetries, but
$M$ admits many non-locally homogenous 
metrics that have all of these symmetries; indeed any $\R^2$-invariant 
perturbation of the $\Sol$-metric has all of these symmetries.  Similar 
examples can be constructed on torus bundles that are very far from 
locally homogenous metrics.

\item Theorem \ref{theorem:general} has an extension to all closed, 
aspherical Riemannian manifolds; in particular, 
normal abelian subgroups are allowed.  See \cite{FW}.
\end{enumerate}

Theorem \ref{theorem:characterize1} can be extended to non-aspherical
manifolds; in the end it is really a theorem about closed manifolds
whose fundamental group is arithmetic.  For a lattice $\Delta$, we 
say that a manifold $M$ has {\em $\Delta$-fundamental hidden symmetries}
if there is a homomorphism $\pi_1(M)\to \Delta$ and hidden symmetries of
$M$ that are intertwined with hidden symmetries of $\Delta$.  When
$\Delta=\pi_1(M)$ then we call these fundamental hidden symmetries.

\begin{theorem}[Hidden symmetry theorem, nonaspherical case]
\label{theorem:nonaspherical}
Let $\Delta$ be a lattice in a semisimple Lie group $G$ with maximal
compact subgroup $K$.  If a closed manifold $M$ has infinitely 
many $\Delta$-fundamental hidden symmetries, then $M$ 
is a Riemannian fibration over the arithmetic manifold
$\Delta\backslash G/K$.
\end{theorem}

In \S\ref{section:counterexample}  
we discuss a ``near counterexample'' to Theorem \ref{theorem:general}
and Theorem
\ref{theorem:nonaspherical}, namely a closed, aspherical 
manifold with infinitely many
``topological hidden symmetries''.  We also state in 
\S\ref{section:counterexample} a conjecture on how metrics
with greater and greater hidden symmetry must degenerate.

It is a pleasure to thank David Fisher and Kevin Whyte for their 
useful comments.

\section{Three proofs of Proposition \ref{proposition:biggest}}
\label{section:3proofs}

Each element $\phi\in \Isom(M,h)$ induces an automorphism of 
$\pi_1(M)$, which (since $\phi$ may move a given basepoint) is only
defined up to conjugacy, thus giving a map $\psi:\Isom(M,h)\to
\Out(\pi_1(M))$.  Mostow Rigidity states that, under the hypothesis of
Proposition \ref{proposition:biggest}, any isomorphism between fundamental
groups of two locally symmetric manifolds is induced by a unique
isometry.  We then have
$$\Isom(M,h)\stackrel{\psi}{\to} \Out(\pi_1M)\rightarrow \Isom(M,g_{loc})$$
with the last map an isomorphism by Mostow Rigidity.  Hence to prove
Proposition \ref{proposition:biggest}, it is enough to prove that $\psi$ is
injective.  

Suppose $\ker(\psi)\neq 0$.  We claim there exists $\phi \in \ker
(\psi)$ with order some prime $p>1$.  To prove this, first note that 
$\Isom(M,h)$ is a compact (since $M$ is compact) Lie group, by the
classical Steenrod-Myers Theorem.  Any element of the connected
component $\Isom^0(M,h)$ 
of the identity of $\Isom(M,h)$ is homotopically trivial, and
so lies in $\ker (\psi)$.  Now every nontrivial, connected, 
compact Lie group contains a circle, and so contains elements of any
prime order.  If $\Isom^0(M,h)$ is trivial, then 
$\Isom(M,h)\supseteq \ker (\psi)$ is a finite group, 
and so contains an element 
$g$ of order some prime $p>1$, and the claim is proved.  

So, we have an element $g\in \Isom(M,h)$ of order $p>1$, with $p$ prime,
so that $g$ is homotopically trivial.  We now give three different 
proofs that $g=1$, proving that $\psi$ must indeed be injective.

\subsection{Smith theory proof}

Let $F=<g>$, and consider the group $\Lambda$ of lifts of every element
of $F$ to $\widetilde{M}$.  There is an exact sequence
\begin{equation}
\label{eq:lift}
1\to \pi_1(M)\to \Lambda\to F\to 1
\end{equation}
with the kernel acting on $\widetilde{M}$ by deck transformations.  
The fact that $g$ is homotopically trivial gives that the action
$F\to \Out(\pi_1(M))$ is trivial.  Further, 
since $M$ is locally symmetric, nonpositively curved, and has 
no local torus factors, it follows that the center $Z(\pi_1(M))=0$, so
in particular $H^2(G,Z(\pi_1(M)))=0$.  These two facts imply that the
exact sequence (\ref{eq:lift}) splits, so that the $F$-action on $M$
lifts to an action by isometries on $(\widetilde{M}, h)$, where by abuse
of notation $h$ denotes the lift to $\widetilde{M}$ of the metric $h$ on
$M$.

Now, the lift of any two metrics on a compact manifold are
quasi-isometric (Effremovich-Milnor-Svarc), so the lift of each element
of $F$ induces a quasi-isometry of $\widetilde{M}$ endowed with the
metric $g_{loc}$.  But each element of $F$ is homotopically trivial on
$M$, and since homotopies are compact, there exists $D>0$ so that 
$$d(x,f(x))\leq D \mbox{\ for all\ }f\in F, x\in \widetilde{M}$$
where $d$ is the distance on $X:=(\widetilde{M}, g_{loc})$.  

Now $X$ has a compactification $\overline{X}=X\cup \partial X$ by
Hausdorff equivalence classes of geodesic rays.  In this
compactification, $\overline{X}$ is homeomorphic to 
a closed $n$-ball, $n=\dim(M)$, and each isometry of
$X$ extends to a homeomorphism of this ball.  Since each $f\in F$ moves
points by at most $D$, and since $\partial X$ is defined by 
Hausdorff equivalence classes of rays, it follows that the homeomorphic
action of any $f\in F$ on $\overline{X}$ restricts to the identity on
$\partial X$.

Now, since $f=g$ has prime order $p$, we may apply 
Smith Theory to the action of $f$ on the homology $n$-ball 
$(\overline{X},\partial X)$, where $n=\dim(M)$. This gives that 
the pair $(\Fix(f),\Fix(f|_{\partial X}))$ is a mod-$p$ homology
$r$-ball for some $0\leq r\leq n$.  Since $f$ is the identity on
$\partial X$, in other words since $\Fix(f|_{\partial X})=\partial X$, it
must be that $r=n$.  Further, if $\Fix(f)\cap X \neq X$, then radial
projection away from a point $z\not \in \Fix(f)\cap X$ would give a homotopy
equivalence of $(\Fix(f),\Fix(f|_{\partial X}))$ with $(\partial X,
\partial X)$, contradicting that the former is a mod-$p$ homology
$n$-disk with $n>1$.  Hence $\Fix(f)\cap X=X$, so that $f$ is trivial on
$X$, and hence on $M$.
\endproof

\subsection{Harmonic maps proof}
\label{subsection:harmonic}

We recall that a map $f:N\to M$ between Riemannian manifolds is {\em
harmonic} if it minimizes the energy functional
$$E(f)=\int_N||Df_x||^2d{\rm vol}_N$$

The key properties of harmonic maps between closed Riemannian 
manifolds which we will need are the following (see, e.g. \cite{SY}, which
is indeed the inspiration for this proof):

\begin{itemize}
\item (Eels-Sampson) 
When the target manifold has nonpositive sectional curvatures, a
harmonic map exists in each homotopy class.

\item (Hartman) The image of any two harmonic mpas in a given homotopy
class are the boundary of an isometric product with an interval; in
particular the harmonic map in a given homotopy class is unique 
when the domain and target have the same dimension.

\item (easy) The precomposition and postcomposition of 
a harmonic map with an isometry gives a harmonic map.

\end{itemize}

Now let $F=<g>$.  Since $(M,g_{loc})$ is nonpositively curved, there is a
unique harmonic map $\phi$ in the homotopy class of the identity map 
$(M,h)\to (M,g_{loc})$.  Precomposing $\phi$ with any element $f\in F$
then gives, since $f$ is homotopically trivial, that 
$$\phi\circ f=\phi$$
so that $\phi$ is $F$-einvariant.  This principal of ``uniqueness
implies invariance'' will be quite useful for us later.

As $\phi$ is $F$-invariant and $F$ is finite, 
the map $\phi$ factors through a map $\psi:M\rightarrow M/F$.  Now, 
restricting to orientation-preserving isometries if necessary, the 
fixed set of $F$ is a submanifold of $M$ of codimension at least two.  
Hence the quotient orbifold $M/F$ is a pseudomanifold, and so one can
apply degree theory to $M/F$.  Since $\phi$ factors through $\psi$, and
since $\psi$ has degree $|F|$, it follows that 
$\phi$ has degree $|F|$.  But $\phi$ is homotopic to the identity
map, and so has degree one.  Hence $|F|=1$, so that $g$ is the
identity.  
\endproof

\subsection{Algebraic-topological proof}

We consider the action of the orbifold  
fundamental group of $M/\Isom(M)$ on the universal cover
$\widetilde{M}$.  Note that this 
group, by (the algebraic form of) Mostow rigidity, has a representation 
into the algebraic group $G$.  
Note that the symmetric space $X=G/K$ is an example of the universal 
space $\underline{EG}$ 
for proper $G$-actions (see e.g. [BCH]), and thus there 
is a well-defined equivariant homotopy class intertwining these
actions.  As in the proof in Subsection \ref{subsection:harmonic}, 
considering the group generated by $\pi_1(M)$ and by $g$ 
gives a contradiction.

\remark
The difference between this proof and the harmonic maps proof is that
in the latter one can pick out a canonical unique map, 
while in the algebraic proof 
there is merely a contractible space of choices.  The proof of our main
result (Theorem \ref{theorem:characterize1}) 
will use the extra rigidity of the analytic proof. 

\section{Proof of the main results}

In this section we prove the main results of this paper.

\subsection{Hidden symmetries in arithmetic manifolds}
\label{subsection:arithmetic:hidden}

We begin this subsection with the following.

\bigskip
\noindent
{\bf Proof of Proposition \ref{proposition:arithmetic}: }
Let $M=\Gamma\backslash G/K$, endowed with the locally symmetric
metric.   The source of the hidden symmetries in the arithmetic case is 
intersections of conjugates of $\Gamma$.  As $M$ is arithmetic, we have that
$[\Comm_G(\Gamma):\Gamma]=\infty$.  Let $g\in
\Comm_G(\Gamma)\setminus N_G(\Gamma)$, where $N_G(\Gamma)$ denotes the
normalizer of $\Gamma$ in $G$ (which is finite, by Bochner).  The 
subgroup $\Delta_g:=g\Gamma
g^{-1}\cap N\Gamma$ then contains a subgroup $H_g$ which is 
normal in $g\Gamma g^{-1}$ but is not normal in $\Gamma$.  The cover of
$M$ corresponding to $H_g$ then has a hidden symmetry.
\endproof

It would be interesting to explore manifolds with hidden symmetries
ocurring in {\em regular} covers.

\bigskip
\noindent
{\bf Hidden symmetries are thin in \boldmath$\Comm_G(\Gamma)$. }
We now give\footnote{We would like to thank Kevin Whyte for this clean 
description.} substance to the claim that the set of hidden symmetries of
$M=\Gamma\backslash G/K$ is a rather thin subset of $\Comm_G(\Gamma)$.

First, note that $g\in G$ is (the lift to $G/K$ of) 
a hidden symmetry of $M$ if and only if $g$ normalizes some finite index
subgroup of $\Gamma$.  Given $g\in \Comm_G(\Gamma)$, if $g$ normalizes 
$\Gamma'$, then some power $g^k$ of $g$ is an inner automorphism of $\Gamma'$ 
since the outer automorphism group of any cocompact lattice in a simple $G\neq
\SL(2,\R)$ is finite by Mostow Rigidity.  Thus multiplying
$g^k$ by an element of $\Gamma'$ we obtain an element which commutes with
$\Gamma'$, and is therefore the identity.  
Thus $g^k\in \Gamma$ for some $k$, which
is certainly not all of $\Comm_G(\Gamma)$.  For an easy example (in the
noncocompact case), most elements of $\SL_n(\Q)$ 
do not have any power in $\SL_n(\Z)$, for example the diagonal 
matrix with entries $(1/2, 2, 1, 1, 1,\ldots)$.

\subsection{Proof of Theorem \ref{theorem:characterize1}}

\bigskip
\noindent
{\bf Non-arithmetic case. }
First suppose that $M=\Gamma\backslash G/K$ 
is not arithmetic.  By a theorem of Margulis, $\Gamma$ has finite index
in $\Delta:=\Comm_G(\Gamma)$.  Since every (the lift to
$\widetilde{M}$) of every hidden symmetry of $M$ lies in $\Delta$, it
follows that $M$ can only have finitely many hidden symmetries.  Now any
perturbation of the metric on $M$ which is invariant with respect to the
finite group $\Delta/\Gamma$ gives a metric on 
$M$ with the same number of hidden symmetries as $(M,g_{loc})$.

\bigskip
\noindent
{\bf Arithmetic case. } Let $X$ denote the universal cover of $(M,h)$.
For any lattice $\Delta$ which intersects $\Gamma=\pi_1(M)$ in a
subgroup of finite index, we may consider the harmonic map $f_\Delta$ in
the homotopy class of the identity map from the cover of $M$
corresponding to $\Delta$ with the metric $h$ to the corresponding cover
with the metric $g_{loc}$.  Uniqueness of harmonic maps gives us, as in
the harmonic maps proof of Proposition \ref{proposition:biggest} in
\S\ref{subsection:harmonic}, that $f_\Delta=f_\Gamma$.

Noting this for all such $\Delta$, we see that the lift $\widetilde{f}$ 
of the harmonic representative $f$ of the identity map $(M,h)\to
(M,g_{loc})$ will be $g$-equivariant with respect any $g$ 
lying in the union of such $\Delta$ in $G$.  Further,
again by uniqueness of harmonic maps, we see that $\widetilde{f}$ will
be equivariant with respect to the group $U$ generated by all such
$\Delta$.  The given implies that $U$ contains $\Gamma$ as a subgroup of
infinite index.  The (topological) closure of $U$ is a closed,
hence Lie, subgroup of $G$.  By the Borel Density theorem (see, e.g., [Ma]),
$G$ is the smallest algebraic subgroup containing $\Gamma$; from
this it follows that $U$ is dense in $G$.

Using the Arzela-Ascoli theorem, for any
$g\in G$, we can find an isometry $\phi$ of $X$ which intrertwines with 
$g$: one simply takes any limit of elements of $U$ converging to
$g$.  We claim that such a $\phi$ is unique.  The only possible
non-uniqueness comes from the stabilizer of a point, a compact
group.  Hence if $\phi$
were not unique one could factor the map $\widetilde{f}:(X,h)\to
(X,g_{loc})$ through some isometry.  If this isometry were finite
order, the order would divide $\deg(f)=1$, a contradiction; if it
were of infinite order, then $f$ would factor through a
lower-dimensional quotient, also contradicting $\deg(f)=1$.  This
proves the claim.

Hence we have shown that $X$ admits an isometric $G$-action, and
that $\widetilde{f}$ intertwines this action with the standard
isometric $G$-action on the symmetric space $(X,g_{loc})$.  Thus, the 
harmonic map $f$ is actually a homothety, and the theorem is proven.
\endproof

As the proof of Theorem \ref{theorem:nonaspherical} is so
similar, we leave it as an exercise to the reader.

\section{Examples}

\subsection{An exotic periodic diffeomorphism}
\label{subsection:exotic}

According to an old example, for any $n>4$ the group $\Z/p\Z$ has a
smooth action on the $n$-sphere $S^n$ with fixed-set a knotted mod-$p$
homology sphere $K$, with $(S^n,K)$ not homeomorphic to $(S^n,S^m)$ for
any $0\leq m\leq n$.  To obtain such an action, take 
a $(k+1)$-dimensional contractible 
manifold $X$ whose product with an interval is a ball, for example Mazur 
manifolds for $k=3$.  Then for $n-k$ even, the manifold 
$X\times D^{n-k}$ has a $\Z/p\Z$ 
action with fixed-point set $X$. 
Restricting to the boundary gives an action an a sphere with fixed-point
set $\partial X$, which is definitely non-simply connected.

Now if one takes a $\Z/p\Z$ action on an $n$-manifold $M$ and if the
fixed-set this action is nonempty, then one can 
equivariantly connect-sum with one of 
the knotted $(S^n,K)$ as above to 
get a new action on $M$ whose fixed set is definitely different, for 
example the fundamental group of the fixed set won't inject into
$\pi_1(M)$.

\subsection{A near counterexample}
\label{section:counterexample}

We now give what seems like a ``near counterexample'' to Theorem
\ref{theorem:general} and Theorem \ref{theorem:nonaspherical}, namely a
closed, aspherical manifold with infinitely many ``topological hidden
symmetries''.

Let $M$ be a noncocompact, arithmetic hyperbolic manifold, say of
dimension $3$, and say with one
end, or ``cusp''.  This end is homeomorphic to the product of a
$2$-torus $T^2$ with $[0,\infty)$.  Let $V$ be the closed manifold 
obtained from $M$ by first deleting $T^2\times (N,\infty)$ for any large
$N$, and then doubling the resulting manifold.  We will not need this,
but Leeb has shown that such $V$ can be given metrics of nonpositive
curvature.  

Since $M$ is arithmetic, it has infinitely many hidden symmetries.  It
is easy to see that any hidden symmetry of $M$ induces a finite order
diffeomorphism of some finite cover of $V$; we call this a ``hidden
topological symmetry''.  Hence $V$ has infinitely many ``topological 
hidden symmetries'', but it is easy to see that $M$ has no locally
symmetric metric, and indeed is not a fiber bundle over any manifold.  
Note that any topological hidden symmetry preserves some Riemannian
metric on the corresponding finite cover (by averaging); of course this
metric depends (as it must, by Theorem \ref{theorem:general}) on the
cover of $M$.

The reason this example ends up not being a counterexample to the
theorems is that the core of a cover is not the cover of the 
core. Equivalently, think about the cover of a cusp.  This will be 
disconnected and different components are likely to cover it different 
numbers of times.  These different covering cusps will be interchanged 
in a hidden symmetry;  this necessitates isogenies of the cusp with 
shifted versions of itself that go further down the tube.

Given this situation, we make the following conjecture.

\begin{conjecture}[Stretching cusps]
\label{conjecture:pinching}
Let $M$ be the double of a cusped arithmetic 
hyperbolic $n$-manifold $N$, or more
generally the double of a negatively curved, arithmetic 
locally symmetric manifold $N$ which is noncompact but has finite volume.  
Suppose $g_i$ is a sequence of Riemannian metrics on $M$ with the
property that $(M,g_i)$ has at least $i$ hidden symmetries.  Then 
the sequence $\{(M,g_i)\}$ Gromov-Hausdorff converges (in the
appropriate sense) to a disjoint union of two isometric copies of $N$.  
\end{conjecture}

Here we are using a notion of convergence where one chooses two
(sequences of) basepoints. It would be interesting to consider the
special case of $N$ hyperbolic and each of the metrics $g_i$ being
nonpositively curved.

\bigskip
\noindent
Dept. of Mathematics, University of Chicago\\
5734 University Ave\\
Chicago, Il 60637\\
E-mail: farb@math.uchicago.edu, shmuel@math.uchicago.edu


\begin{thebibliography}{ABCDEF}
\small

\bibitem[E]{E}
P. Eberlein, Lattices in spaces of nonpositive curvature, {\em Annals of 
Math.} (2) 111  (1980), no. 3, 435--476.

\bibitem[Fu]{Fu}
A. Furman, Mostow-Margulis rigidity with locally compact targets, 
{\em Geom. Funct. Anal.} 11 (2001), no. 1, 30--59.

\bibitem[FW]{FW}
B. Farb and S. Weinberger, Hidden symmetry, in 
preparation.

\bibitem[Ma]{Ma}
G.A. Margulis, {\em Discrete subgroups of semisimple Lie groups}, Ergebnisse 
der Mathematik und ihrer Grenzgebiete (3), 17. Springer-Verlag, Berlin,
1991.

\bibitem[SY]{SY} 
R.Schoen and S.T.Yau, Compact group actions and the topology of 
manifolds with nonpositive curvature, {\em Topology}  18  (1979), no. 4, 
361--380.

\end{thebibliography}
\end{document}